\documentclass[12pt,a4paper]{article}

\usepackage{amsfonts}
\usepackage{amsmath}
\usepackage{amssymb}
\usepackage{amstext}
\usepackage{amsthm}   

\newcommand{\R}{\mathbb R}
\newcommand{\N}{\mathbb N}

\newcommand{\W}{{\cal{W}}}
\newcommand{\V}{{\cal{V}}}
\newcommand{\C}{{\cal{C}}}

\newtheorem{teo}{Theorem}

\newtheorem{lemma}[teo]{Lemma}

\newtheorem{defin}{Definition}

\newtheorem*{lem}{Lemma A.1}

\theoremstyle{definition}

\theoremstyle{remark}
\newtheorem{rem}[teo]{Remark}

\DeclareMathOperator{\codim}{codim}
 
\DeclareMathOperator{\mis}{meas}

\title{A multiplicity result
for the problem $\delta d \xi = f'(\langle\xi,\xi\rangle )\xi$}
\author{Antonio Azzollini}
\date{Dipartimento di Matematica\\
Universit\`a di Bari,\\
Via Orabona 4, 70125 Bari, Italy\\
e-mail: azzollini@dm.uniba.it}

\begin{document}
\maketitle
\begin{abstract}
In this paper we consider the nonlinear equation involving
differential forms on a compact Riemannian manifold $\delta d \xi
= f'(\langle\xi,\xi\rangle )\xi$. This equation is a
generalization of the semilinear Maxwell equations recently
introduced in a paper by Benci and Fortunato. We obtain a
multiplicity result both in the positive mass case (i.e.
$f'(t)\geq\varepsilon>0$ uniformly) and in the zero mass case
($f'(t)\geq 0$ and $f'(0)=0$) where a strong convexity hypothesis
on the nonlinearity is assumed.
\end{abstract}

\section*{Introduction}
Let $(M,g)$ be a compact Riemannian $n-$manifold, where $n\geq 3$,
and $\Lambda^k(M)$ be the set of regular
 $k$-forms on $M$. We consider the following equation
\begin{equation}\label{problem}
\left\{
\begin{array}{ll}
\delta d\xi=f'(\langle\xi,\xi\rangle)\xi,\\
\\
\xi\in\Lambda^k(M),&1\leq k\leq n-1,
\end{array}
\right.
\end{equation}
where $f:\R\rightarrow\R$ is a $C^2$ map, $d$ is the exterior
differential, $\delta$ is its adjoint with respect to the inner
product
\begin{equation}\label{scalarproduct}
(\eta
,\xi)_2=\int_M\langle\eta,\xi\rangle\,\omega=\int_M*(\eta\wedge*\xi)\,\omega
\end{equation}
$*$ is the Hodge operator and $\omega$ is a volume n-form.
\\
In this paper we are looking for weak solutions of
\eqref{problem}, namely for solutions of
\begin{equation}\label{weakproblem}
\left\{
\begin{array}{l}
\xi\in H^1_k(M),\hbox{(see section \ref{property} for the
definition)}\\
\\
\int_M\langle d\xi,d\eta\rangle\,\omega=\int_M
f'(\langle\xi,\xi\rangle)\langle\xi,\eta\rangle\,\omega,\quad\forall\eta\in
H^1_k(M).
\end{array} \right.
\end{equation}
\newline
If we set
    \begin{equation}\label{defF}
        F(\xi):=\int_Mf(\langle\xi,\xi\rangle)\omega,
    \end{equation}
then, assuming a suitable condition on the growth of $f'$, by
standard arguments we have that $F\in C^1(H_k^1(M)),$
  so in order to solve
\eqref{weakproblem} we find critical points of the functional
\begin{equation}\label{J}
J(\xi)=\int_M\langle d\xi,d\xi\rangle\,\omega-F(\xi)
\end{equation}
defined for all $\xi\in H^1_k(M)$.

The strongly indefinite nature of the functional $J$, largely
discussed in \cite{B.F.}, doesn't allow us to approach this
problem in a standard way. In other words, the functional $J$
doesn't present the geometry
of the mountain pass in any space with finite codimension.\\
Assume that
    \begin{description}
        \item{$f_1$)} $f(0)=0$, and $\exists\, \varepsilon >0$ s.t.
        $\forall t \geq 0: f'(t)\geq
        \varepsilon$,
        \item{$f_2$)} $f$ is strictly convex,
    \end{description}
and for $p\in ]2,\frac{2n}{n-2}[$
    \begin{description}
        \item{$f_3$)} $\exists a>0, b>0$ s.t. $|f'(t)|\leq a\,
        t^{\frac{p}{2}-1}+b$, $\forall t\geq 0$,
        \item{$f_4$)} $\exists R>0$ s.t.
        $0<\frac{p}{2}\,f(t)\leq f'(t)t$ for $t>R.$
    \end{description}
We have the following result
\begin{teo}\label{main}
If $f_1$)-$\ldots$-$f_4$) hold, then the problem
\eqref{weakproblem} has infinitely many solutions.

Moreover the same conclusion holds if $f_1$) and $f_2$) are
substituted respectively by
    \begin{description}
        \item{$\tilde f_1$)} $f(0)=f'(0)=0$, and $\forall t \geq 0: f'(t)\geq
        0$,
        \item{$\tilde f_2$)} $\exists \overline c>0$ s.t. $\forall \xi,\eta \in
        \Lambda^k(M)$
            $$f(\langle \xi,\xi\rangle)-f(\langle \eta,\eta\rangle)-
            2 f'(\langle
            \eta,\eta\rangle)\langle\eta,\xi-\eta\rangle\geq\overline
            c\langle \xi-\eta,\xi-\eta\rangle^\frac{p}{2}$$
        pointwise in $M.$
    \end{description}
\end{teo}
\begin{rem}
As in \cite{B.L.}, in the sequel we shall refer to the hypotheses
$f_1$ and $\tilde f_1$ respectively as the ``positive mass" and
``zero mass" case.

Moreover we want to point out the fact that $\tilde f_2$ is just a
pointwise convexity condition. In fact for every $q\in M$ we can
define the scalar product $\langle\cdot,\cdot\rangle_q$ on the
vector space $\Lambda^k(M)$ and the functional
    $$I_q(\xi)=f(\langle\xi,\xi\rangle_q).$$
Since $I_q'(\eta)=2f'(\langle
            \eta,\eta\rangle)\langle\eta,\cdot\rangle$, $\tilde f_2$ implies that
for all $\xi,\eta \in
        \Lambda^k(M)$ s.t. $\eta\neq\xi$
        $$I_q(\xi)-I_q(\eta)-\langle
        I_q'(\eta),\xi-\eta\rangle>0,$$
and then $I_q$ is strictly convex. In the Appendix we shall show
that $\tilde f_2$ is satisfied when $f(t)=t^\frac{p}{2}$.
\end{rem}
\bigskip In order to prove Theorem \ref{main}, as in \cite{B.F.} we shall use
Hodge decomposition to split $\xi$ in this way
\begin{equation}
\xi=d\alpha+\beta
\end{equation}
where $\delta\beta=0$.

The functional \eqref{J} formally becomes
    \begin{equation}\label{J2}
        J(\alpha,\beta) = \int_M\langle
d\beta,d\beta\rangle\,\omega- \int_M f(\langle
d\alpha+\beta,d\alpha+\beta\rangle)\,\omega.
\end{equation}
If we set
    \begin{eqnarray}
        J_\alpha&:& \beta\mapsto J_\alpha(\beta)=J(\alpha,\beta)\\
        \noalign{\hbox{and}}\nonumber\\
        J_\beta&:&\alpha\mapsto J_\beta(\alpha)=J(\alpha,\beta),
    \end{eqnarray}
then we can define the partial derivative of $J$ as follows
    \begin{eqnarray}
        \frac{\partial J}{\partial\alpha}(\alpha,\beta): & = &
        dJ_\beta(\alpha)\\
        \frac{\partial J}{\partial\beta}(\alpha,\beta): & = &
        dJ_\alpha(\beta).
    \end{eqnarray}
We are interested in finding critical points of \eqref{J2}, i.e.
the couples $(\alpha,\beta)$ such that
\begin{eqnarray}\label{partialalpha}
\frac{\partial J}{\partial\alpha}(\alpha,\beta) & = & 0\\
\noalign{\hbox{and}}\nonumber\\ \frac{\partial
J}{\partial\beta}(\alpha,\beta) & = & 0.\label{partialbeta}
\end{eqnarray}
Set
\begin{equation}\label{Fbeta}
F_\beta:  =  \alpha\mapsto F(d\alpha+\beta)
\end{equation}
and note that by the convex nature of \eqref{Fbeta}, the problem
\eqref{partialalpha} is actually a minimizing problem.

In section \ref{first part}, where we assume $f_1)\ldots f_4)$, we
introduce some preliminary results and the definition of the
spaces $\V$ and $\W$ which respectively $\alpha$ and $\beta$
belong to. These spaces are constructed in such a way we have, for
any $\beta\in\W$, a unique solution $\Phi(\beta)\in\V$ for the
minimizing problem
\eqref{partialalpha}. \\
Then in Theorem \ref{critpoin} we'll show that in order to solve
the system \eqref{partialalpha} and \eqref{partialbeta} we are
reduced to study the critical points of the functional
    \begin{equation}\label{totalfunc}
        \widehat{J}(\beta)=J(\Phi(\beta),\beta),\qquad\beta\in\W.
    \end{equation}
Differently from $J$, $\widehat J$ doesn't exhibit strong
indefinitness, so the proof of the existence of infinitely many
critical points is carried out by using a well known multiplicity
result for even functionals.

In section \ref{zeromass} we replace $f_1$ and $f_2$ by $\tilde
f_1$ and $\tilde f_2$. Differently from the previous situation, we
can't give any proof about the regularity of the functional
\eqref{totalfunc}. To overcome this difficulty we work in an
indirect way. In fact we perturb the problem \eqref{problem} by
adding a linear ``mass term" $m\xi$, $m>0$, to the nonlinearity on
the right hand side. For this perturbed problem we have infinitely
many solutions since it satisfies the
hypothesis $f_1$.\\
Then we study the behaviour of the solutions of the perturbed
problem when the perturbation goes to zero.

\section{Positive mass case}\label{first part}

\subsection{The functional framework}
For any $q>1$ and $k\in \N$, let $H^{1,q}_k(M),$ $H^1_k(M)$ and
$L^q_k(M)$ be defined as follows
\begin{eqnarray*}
H^{1,q}_k(M)    &:=&\overline{\Lambda^k(M)}^{\|\cdot\|_{1,q}}\;,\\
H^1_k(M)        &:=&\overline{\Lambda^k(M)}^{\|\cdot\|},\\
L^q_k(M)        &:=&\overline{\Lambda^k(M)}^{|\cdot|_q}
\end{eqnarray*}
where, for every $\xi\in\Lambda^k(M)$,
\begin{eqnarray}
\|\xi\|^q_{1,q} &:=& \int_M\langle d \xi,d
\xi\rangle^{q/2}\,\omega+\int_M\langle\delta\xi,\delta\xi\rangle^{q/2}\,\omega+
\int_M\langle\xi,\xi\rangle^{q/2}\,\omega\nonumber\\
\|\xi\|^2       &:=& \int_M\langle d \xi,d
\xi\rangle\,\omega+\int_M\langle\delta\xi,\delta\xi\rangle\,\omega+
\int_M\langle\xi,\xi\rangle\,\omega\label{normH}\\
|\xi|_q^q &:=&\int_M\langle\xi,\xi\rangle^{q/2}\,\omega.\nonumber
\end{eqnarray}
Since
    \begin{equation}\label{immer}
        H^1_k(M)\hookrightarrow L^q_k(M)\quad \hbox{for }
        1\leq q\leq\frac{2n}{n-2},
    \end{equation}
by $f_3$ we have that $J(\alpha,\beta)<+\infty$ for $\alpha\in
H^{1,p}_{k-1}(M)$ and $\beta\in H^1_k(M)$.\\

Observe that $F_\beta:H_{k-1}^{1,p}(M)\rightarrow\R$ defined by
\eqref{Fbeta} is not coercive, since it is constant on the space
    \begin{equation}\label{C}
        \C:=\left\{\eta\in H^{1,p}_{k-1}(M)\mid d\eta=0\right\}.
    \end{equation}
However we have the following result
    \begin{lemma}\label{coerciveness}
        For all $\beta\in L^p_k(M),$  $F_\beta$ is coercive on
        \begin{equation}
            \V:=\left\{\alpha\in H^{1,p}_{k-1}(M)\mid\forall\eta\in\C:\int_M\langle\alpha,
            \eta\rangle\,\omega=0\right\}.
\end{equation}

\begin{proof}
First observe that, by $f_4,$
    \begin{equation}\label{increasing}
        \exists c>0, d>0  \hbox{ s.t. }
        c\,t^{\frac{p}{2}}\leq f(t)+d,  \forall t\geq 0.
    \end{equation}
Let $\beta\in L^p_k(M)$. By Lemma $6$ in \cite{A}, we know that
the norm on $H_k^1(M)$ defined by
    \begin{equation}\label{normnorm}
        \|\xi\|_\sim^2:=\int_M(\langle d\xi,d
        \xi\rangle+\langle\delta\xi,\delta\xi\rangle+\langle\xi^0,\xi^0\rangle)\,\omega
    \end{equation}
where $\xi^0$ is the orthogonal projection of $\xi$ on
$\ker(-\Delta)$, is equivalent to the norm defined by
\eqref{normH}.\\
In particular, in the space $\V$ the norm \eqref{normnorm} becomes
    \begin{equation*}
        \|\xi\|_\sim^2=\int_M\langle d\xi,d
        \xi\rangle\,\omega.
    \end{equation*}
Indeed, if $\alpha\in\V,$ then
    \begin{equation}\label{deltaalpha=0}
        \int_M\langle\delta\alpha,\delta\alpha\rangle\,\omega=
        \int_M\langle d\delta\alpha,\alpha\rangle\,\omega=0
    \end{equation}
because $d\delta\alpha\in\C.$\\
Moreover, since $\alpha^0\in\ker(-\Delta),$ then $\alpha^0\in\C.$
But $\alpha^0\in\V$ and then
    \begin{equation}\label{alpha0=0}
        \int_M\langle\alpha^0,\alpha^0\rangle\,\omega=0.
    \end{equation}
By \eqref{deltaalpha=0} and \eqref{alpha0=0} we can conclude that
$|d\alpha|_2$ is a norm on the space $\V$ equivalent to
$\|\alpha\|$, i.e.
    \begin{equation}\label{eqnorma}
        \exists \tilde c>0\; s.t.\;\forall\alpha\in\V:\|\alpha\|\leq\tilde
        c|d\alpha|_2.
    \end{equation}
By \eqref{eqnorma} and since $H^1_{k-1}(M)\hookrightarrow
L^p_{k-1}(M)$,
$$\|\alpha\|^p_{1,p}=|d\alpha|_p^p+|\alpha|_p^p\leq|d\alpha|_p^p+
c_1\|\alpha\|^p\leq|d\alpha|_p^p+c_2|d\alpha|_2^p$$
and then
    \begin{equation}\label{bo}
        \|\alpha\|_{1,p}\leq c_3(|d\alpha|_p+|d\alpha|_2)\leq c_4|d\alpha|_p
    \end{equation}
since $L_k^p(M)\hookrightarrow L_k^2(M).$\\
Now, consider $(\alpha_n)_{n\geq1}$ in $\V$ s.t.
$\|\alpha_n\|_{1,p}\rightarrow+\infty.$ By \eqref{bo}
    \begin{equation}\label{coer1}
        \int_M\langle d\alpha_n+\beta,d\alpha_n+\beta\rangle^\frac{p}{2}\,\omega
        \rightarrow+\infty,
    \end{equation}
so by \eqref{increasing} we have
    \begin{equation}\label{coer2}
        c\int_M\langle d\alpha_n+\beta,d\alpha_n+\beta\rangle^\frac{p}{2}\,\omega
        \leq d \mis(M)+\int_Mf(\langle
        d\alpha_n+\beta,d\alpha_n+\beta\rangle)\,\omega.
    \end{equation}
The coerciveness of $F_\beta$ is a consequence of \eqref{coer1}
and \eqref{coer2}.
\end{proof}
\end{lemma}

\begin{teo}
For every $\beta\in L^p_k(M)$ there exists a unique minimizer of
$F_{\beta|\V}$.
\end{teo}
\begin{proof}
    Let $\beta\in L^p_k(M)$. By $f_2$ and $f_3$ the functional
    $F:L^p_k(M)\rightarrow\R$ defined by \eqref{defF} is strictly
    convex and continuous. Obviously, also $F_\beta$ has the same
    properties, so it is weakly lower semicontinuos. Since
    $F_\beta$ is also coercive in $\V$ by Lemma \ref{coerciveness}, certainly it
    possesses a minimizer $\Phi(\beta)\in\V$. The uniqueness is a consequence of the
    strict convexity.
\end{proof}

\smallskip
\subsection{Regularity, symmetry and compactness}\label{property}
Assume the following definitions:
\begin{equation}\label{Phi}
\begin{array}{cll}
\Phi: L^p_k(M)\rightarrow\V & s.t. &\Phi(\beta)\hbox{ is the
minimizer of }F_{\beta|\V}\\
\\
\widehat{J}:\W\rightarrow\R & s.t.
&\forall\beta\in\W:\widehat{J}(\beta)=J(\Phi(\beta),\beta)
\end{array}
\end{equation}
where $$\W:=\left\{\beta\in H^1_k(M)\mid \delta\beta=0\right\}.$$
\begin{teo}\label{critpoin}
If $\Phi\in C^1(\W,\V)$, then $\widehat{J}\in C^1(\W)$ and its
critical points  are solutions of \eqref{partialalpha},
\eqref{partialbeta}.
\end{teo}
\begin{proof}
Suppose $\Phi\in C^1(\W,\V),$ then certainly $\widehat{J}\in
C^1(\W)$ since it is the composition of $C^1$ maps.\\
Now let $\beta_0\in\W$ be a critical point of $\widehat{J}$. We
have that for any $\beta\in\W$
\begin{equation}
0=\langle\widehat{J}\,'(\beta_0),\beta\rangle=\big\langle\frac{\partial
J}{\partial\alpha}(\Phi(\beta_0),\beta_0),\Phi'(\beta_0)(\beta)\big\rangle+
\big\langle\frac{\partial
J}{\partial\beta}(\Phi(\beta_0),\beta_0),\beta\big\rangle,
\end{equation}
that is
\begin{equation}
\frac{\partial
J}{\partial\beta}(\Phi(\beta_0),\beta_0)=-\frac{\partial
J}{\partial\alpha}(\Phi(\beta_0),\beta_0)\circ\Phi'(\beta_0).
\end{equation}
But
\begin{equation}\label{partialalpha=0}
\frac{\partial J}{\partial\alpha}(\Phi(\beta_0),\beta_0)=0
\end{equation}
because $\Phi(\beta_0)$ is a minimizer of $F_{\beta|\V}$, so also
\begin{equation}
\frac{\partial J}{\partial\beta}(\Phi(\beta_0),\beta_0)=0.
\end{equation}
\end{proof}
In order to study the functional $\widehat J$, we need to
investigate the properties of the map $\Phi$. Then, in the next
theorem we are going to prove some regularity, symmetry and
compactness properties of the map $\Phi$. To get regularity, in
particular, we will use the implicit function theorem on
\begin{equation}
    \frac{\partial J}{\partial\alpha}:\V\times\W\rightarrow\V'.
\end{equation}
where $\V'$ is the dual of $\V$.\\
Observe that $J\in C^2(\V\times\W)$. Moreover we have the
following
\begin{lemma}\label{uniformconvex}
Set
        \begin{equation}
            \tilde F:=(\alpha,\beta)\in\V\times\W\mapsto
            F(d\alpha+\beta).
        \end{equation}
Then for all $\beta\in L^p_k(M)$, $\tilde F(\cdot,\beta)\in
C^2(\V)$ is uniformly convex on $\V$ with respect to the norm
$\|\cdot\|$, i.e. there exists $C>0$ such that for all $\alpha$,
$\overline\alpha\in\V$
    \begin{equation}
        \frac{\partial^2 \tilde F}{\partial\alpha^2}(\alpha,\beta)
        [\overline\alpha,\overline\alpha]:
        =d^2 \tilde F(\cdot,\beta)(\alpha)
        [\overline\alpha,\overline\alpha]
        \geq C\|\overline\alpha\|^2.
    \end{equation}
\begin{proof}
Let $\beta\in L^p_k(M)$ and $\alpha\in\V.$ If
$\overline\alpha\in\V,$ then
    \begin{eqnarray*}
        \frac{\partial^2 \tilde F}{\partial\alpha^2}(\alpha,\beta)[\overline\alpha,\overline\alpha]
        &=& 4\int_M
        f''(\langle d\alpha+\beta,d\alpha+\beta\rangle)(\langle d\alpha+
        \beta,d\overline\alpha\rangle)^2\,\omega\\
        &&\phantom{'()^2\beta}+2\int_M f'(\langle
        d\alpha+\beta,d\alpha+\beta\rangle)\langle
        d\overline\alpha,d\overline\alpha\rangle\,\omega\\
        &\geq&\varepsilon\int_M\langle
d\overline\alpha,d\overline\alpha\rangle\,\omega,
\end{eqnarray*}
where the last inequality follows from the convexity of $f$ and
assumption $f_1.$
\end{proof}
\end{lemma}
Now we are ready to prove
\begin{teo}\label{propertyphi}
The following properties hold

\begin{enumerate}
    \item $\Phi\in C^1(\W)$;
    \item $\Phi$ is odd;
    \item $\Phi$ is compact.
\end{enumerate}
\begin{proof}
\begin{enumerate}
    \item
        Since $\Phi(\beta)$ is a minimizer of
        $F_\beta$
        \begin{equation}\label{partderivJ=0inPhibeta}
            \forall\beta\in\W:\frac{\partial
            \tilde F}{\partial\alpha}(\Phi(\beta),\beta)=0,
        \end{equation}
        so, by Lemma \ref{uniformconvex} and the implicit function theorem,
        we have that $\Phi\in C^1(\W,\V)$ and for any
        $\beta\in\W$:
        $$\Phi'(\beta)=-\left(\frac{\partial^2 \tilde F}{(\partial\alpha)^2}
        (\Phi(\beta),\beta)\right)^{-1}\circ
        \frac{\partial^2 \tilde F}{\partial\alpha\partial\beta}(\Phi(\beta),\beta).$$
    \item If $\beta\in\W,$ then some calculations show that
        $$\frac{\partial
        \tilde F}{\partial\alpha}(-\Phi(\beta),-\beta)=\frac{\partial
        \tilde F}{\partial\alpha}(\Phi(\beta),\beta)=0,$$
        so, by uniqueness, $-\Phi(\beta)=\Phi(-\beta).$
    \item  By the same arguments as above, we can prove that $\Phi$ is also in
    $C^1(L^p_k(M),\V)$ so,
        if $(\beta_n)_n$ is a sequence in $\W$ and
        $$\beta_n\rightharpoonup\beta\in\W$$
        with respect to the norm $\|\cdot\|,$ then
        $$\beta_n\rightarrow\beta \hbox{ in }L^p_k(M)$$
        and so
        $$ \Phi(\beta_n)\rightarrow\Phi(\beta)\hbox{ in }\V.$$
\end{enumerate}
\end{proof}

\end{teo}

\subsection{Main Theorem (first part)}\label{proof}
We introduce some results on the Laplace Beltrami operator
$-\Delta$.

It is well known that $-\Delta$ is a self adjoint operator on $L^2_k(M)$ with a
nonnegative, discrete and divergent spectrum
$\sigma(\Delta)$.

Set
    \begin{equation}
    \lambda_1\leq\lambda_2\leq\ldots
    \end{equation}
the sequence of the eigenvalues different from zero repeated
according to their finite multiplicity. The corresponding
eigenvectors
    \begin{equation}
    \eta_1,\eta_2,\ldots
    \end{equation}
constitute an orthonormal basis of $(\ker(-\Delta))^\perp$. Take a
basis $h_1,h_2,\ldots,h_N$ of $\ker(-\Delta)$ so that, if
$\beta\in L^2_k(M)$, we have
    $$\beta=\sum_{i\geq 1}^\infty\beta_i \eta_i+\beta^0,$$
where $(\beta_i)_{i\geq 1}$ are the Fourier components of $\beta$
corresponding to $\eta_1,\eta_2\ldots$, and $\beta^0$ is its
projection on $\ker(\Delta)$.

For every $s\in\R$ define the Sobolev space
    \begin{eqnarray}\label{Bessel2}
        W^{s,2}_k(M)&:=&\Big\{\beta\in L^2_k(M)\mid
        \|\beta\|_{s,2}<\infty\Big\},\nonumber\\
\noalign{\hbox{where}}
        \|\beta\|^2_{s,2}&:=&\sum_{i\geq
        1}^\infty (\lambda_i^s+1)\beta_i^2+|\beta^0|_2^2.
    \end{eqnarray}
It is easy to prove that $(W^{s,2}_k(M),\|\cdot\|_{s,2})$ is a
Banach space and $W^{1,2}_k(M)\equiv H^1_k(M).$ Since
$p<\frac{2n}{n-2}$, by Sobolev embedding theorem, there exist
$s<1$ and $\tilde c>0$ such that
    \begin{equation}\label{Sobolev}
        |\beta|^2_p\leq \tilde c\|\beta\|^2_{s,2},\qquad\forall\beta\in W_k^{s,2}(M).
    \end{equation}

Now we recall the following abstract multiplicity theorem whose
proof can be found in \cite{B.B.F.} (see also \cite{A.R.}).
\begin{teo}\label{mult}
Let $H$ be an Hilbert space and $I$ be a $C^1$ even functional on
$H$ such that
\begin{enumerate}
    \item I(0)=0,
    \item I satisfies (P-S) condition i.e. any sequence $(x_n)_n$ such that
        $$
        \begin{array}{c}
            I(x_n) \hbox{ is bounded}\\
            \\
            I'(x_n)\rightarrow0,
        \end{array}
        $$
        admits a convergent subsequence,
    \item there exist $H^-$, $H^+$ two closed subspaces of $H$ such
    that
    \begin{enumerate}
        \item $\codim (H^+)<\dim (H^-)<+\infty$
        \item $\exists\, c_0,\rho>0$ s.t. $ I(x)\geq c_0,$ $\forall
        x\in\partial S_\rho(0)\cap H^+$, $($where $\partial S_\rho(0):=
        \left\{x\in H\mid\|x\|=\rho\right\})$
        \item $\exists\, c_1>0$ s.t. $\forall x\in H^-:I(x)<c_1$
    \end{enumerate}
\end{enumerate}
then $I$ possesses at least $\dim (H^-)- \codim (H^+)$ couples of
critical points whose corresponding critical values are in
$[c_0,c_1].$
\end{teo}

In the next lemmas we shall verify the hypotheses of the previous
theorem for the functional $\widehat J$ on the Hilbert space $\W$.

    \begin{lemma}\label{ps condition}
        $\widehat J$ is a $C^1$ even functional satisfying the
        (P.S.) condition
    \end{lemma}
    \begin{proof}
        The regularity and symmetry properties of $\widehat{J}$ are an
            immediate consequence of the structure of $J,$ and Theorem \ref{propertyphi}.

        As to (P.S.) condition, let $(\beta_n)_n$ be a sequence in $\W$ such that
            \begin{eqnarray}
                \widehat{J}(\beta_n) & = & |d\beta_n|_2^2-\int_M f(b_n)\,\omega\leq
                M\label{bound},\;M\geq0\\
                \noalign{\hbox{and}}\nonumber\\
                \widehat J\,'(\beta_n) & \longrightarrow & 0\label{diffgo0}
            \end{eqnarray}
        where we have set $b_n=\langle\beta_n+d\Phi(\beta_n),\beta_n+d\Phi(\beta_n)\rangle$
        to simplify the notations.
        We want to show that $\{\beta_n\}$ is precompact. By
        using $f_3$ and $3$ of Theorem \ref{propertyphi} it can be easily seen that $\widehat
        J\,'$ is the sum of an homeomorphism and a compact map, so, by
        standard arguments, we are reduced to prove that
        $(\beta_n)_n$ is bounded.

        If we set
            \begin{eqnarray}
                \varepsilon_n:&=&\frac{1}{2}\Big\langle\widehat{J}\,'(\beta_n),
                \frac{\beta_n}{\|\beta_n\|}\Big\rangle,\label{defvareps}\\
                \noalign{\hbox{from \eqref{diffgo0} we deduce that}}\nonumber
                \varepsilon_n&\longrightarrow&0.\label{diffgo02}
            \end{eqnarray}

        Rendering \eqref{defvareps} explicit, we obtain
            \begin{equation}\label{convergence}
                |d\beta_n|_2^2-\int_M f'(b_n)\langle\beta_n+d\Phi(\beta_n),
                \beta_n+d\big(\Phi'(\beta_n)(\beta_n)\big)\rangle\,\omega=
                \varepsilon_n\|\beta_n\|.
            \end{equation}

    By \eqref{partderivJ=0inPhibeta} we have
        \begin{eqnarray}
            \int_M
            f'(b_n)\langle\beta_n+d\Phi(\beta_n),d\big(\Phi'(\beta_n)(\beta_n)\big)\rangle\,
            \omega & = &
            0\label{ping}
            \\
            \int_M
            f'(b_n)\langle\beta_n+d\Phi(\beta_n),d\Phi(\beta_n)\rangle\,\omega
            & = & 0\label{pong},
        \end{eqnarray}
    so, comparing \eqref{convergence}, \eqref{ping} and \eqref{pong} we
    have
        \begin{equation}\label{convergence2}
            |d\beta_n|_2^2
            -\int_M f'(b_n)\langle\beta_n+d\Phi(\beta_n),
            \beta_n+d\Phi(\beta_n)\rangle\,\omega
            =\varepsilon_n\|\beta_n\|.
        \end{equation}
        Comparing \eqref{bound} and \eqref{convergence2} we get,
    \begin{equation}\label{first}
        \frac{p-2}{2}|d\beta_n|_2^2-
        \int_M \Big(\frac{p}{2} f(b_n) - f'(b_n)b_n\Big)\,\omega \leq
        K_1+K_2\|\beta_n\|
    \end{equation}
and
    \begin{equation}\label{second}
       \int_M \Big(f'(b_n)b_n - f(b_n)  \Big)\,\omega\leq
        K_3+K_4\|\beta_n\|
    \end{equation}
where $K_1,\ldots,K_4$ and the following $\left\{(K_i)\mid
i\geq5\right\}$ are positive constants.\\
By $f_4$, there exists $L>0$ such that
    \begin{equation}\label{second2}
        -L\leq f'(t)t-\frac{p}{2}f(t),\quad\forall t\geq 0,
    \end{equation}
so, by \eqref{increasing},
    \begin{equation}\label{second3}
        f'(t)t-f(t)\geq -L+\frac{p-2}{2}f(t)\geq-K_5+
        K_6\,t^{\frac{p}{2}},\quad\forall t\geq 0.
    \end{equation}
From \eqref{second3} and \eqref{second} we have
    \begin{equation}\label{fourth}
        |\beta_n+d\Phi(\beta_n)|_p\leq
        K_7+K_8\|\beta_n\|^\frac{1}{p},
    \end{equation}
so, since $L^p_k(M)\hookrightarrow L^2_k(M),$
    \begin{eqnarray}\label{fourth2}
        |\beta_n|_2^2 & \leq &
        |\beta_n|_2^2+|d\Phi(\beta_n)|_2^2=|\beta_n+d\Phi(\beta_n)|_2^2\nonumber\\
                        & \leq & K_9|\beta_n+d\Phi(\beta_n)|_p^2\leq
        K_{10}+K_{11}\|\beta_n\|^\frac{2}{p}.
    \end{eqnarray}
From \eqref{first} and \eqref{second2} we also derive
    \begin{equation}\label{third}
        |d\beta_n|_2^2\leq
        K_{12}+K_{13}\|\beta_n\|.
    \end{equation}
Inequalities \eqref{fourth2} and \eqref{third} imply that the
sequence $(\|\beta_n\|)_n$ is bounded.
    \end{proof}
\smallskip
Now, for any $\mu>0$ and $\rho>0$, we set
    \begin{eqnarray}
        \partial S_\rho & := &
        \left\{\beta\in\W\mid\|\beta\|=\rho\right\},\\
        H^+(\mu) & := &
        \overline{\bigoplus_{\lambda_i>\mu}M_{\lambda_i}},\label{H+}\\
        H^-(\mu) & := & (H{^+}(\mu))^\perp\oplus
        M_{\lambda_k},\label{H-}
    \end{eqnarray}
where $(M_{\lambda_i})_{i\geq 1}$ are the spaces of the
eigenfunctions corresponding to the eigenvalues
$(\lambda_i)_{i\geq 1}$ and $k:=\min\left\{i\in \N\mid
\lambda_i>\mu\right\}.$ Observe that, since every eigenspace has
finite dimension, from \eqref{H+} and \eqref{H-} we deduce
    \begin{equation}\label{3a}
        \dim H^-(\mu)=\codim H^+(\mu)+\dim M_{\lambda_k}<+\infty.
    \end{equation}
Moreover we have that
    \begin{lemma}\label{geometry}
        There exist a strictly increasing sequence $(C_i)_{i\geq1}$ and
        two positive numbers sequences $(\rho_i)_{i\geq1}$ and
        $(\mu_i)_{i\geq1}$ such that for
        all $i\geq1$ we have
        \begin{description}
            \item{3b)} $ \widehat J(\beta)\geq C_{2i},$ $\forall
            \beta\in\partial S_{\rho_i}\cap H^+(\mu_i)$,
            \item{3c)} $\sup_{\beta\in H^-(\mu_i)}\widehat
            J(\beta)<C_{2i+1}$.
        \end{description}
    \end{lemma}
    \begin{proof}
        We will prove that for every $C>0$ there exist $\mu>0$
        and $\rho>0$ such that
            \begin{eqnarray}
                \widehat J(\beta)&\geq& C,\quad\forall
                \beta\in\partial S_\rho\cap H^+(\mu)\label{geometry1}\\
                \sup_{\beta\in H^-(\mu)}\widehat J(\beta) & < &
                +\infty.\label{geometry2}
            \end{eqnarray}
        Set $C>0$.   \\
        By $f_3$ and \eqref{Phi}, we have that there exists $b'>0$ such
        that
    \begin{eqnarray}\label{non lin bound}
        \int_M f(\langle \beta+d\Phi(\beta),\beta+d\Phi(\beta)\rangle)\,\omega
        & \leq & \int_Mf(\langle\beta,\beta\rangle)\,\omega\nonumber\\
        & \leq & a |\beta|_p^p+b|\beta|_2^2\nonumber\\
        & \leq & a |\beta|_p^p+b'|\beta|_p^2.
    \end{eqnarray}
    Set
    $\mu=\rho^{\frac{2}{1-s}}$, and
    $K=\displaystyle{\min_{\lambda_i>\mu}}$$\frac{\lambda_i^s}{\lambda_i^s+1}>0$ where
    $\rho$ is a suitable real number
    that we are going to evaluate and $s\in(0,1)$ is defined as in \eqref{Sobolev}.

    Let $\beta\in H^+(\mu)\cap S_\rho$. Using \eqref{Bessel2}, \eqref{Sobolev}
    and \eqref{immer},   we have
        \begin{eqnarray}\label{control}
           |d\beta|^2_2 & = &
          \sum_{\lambda_i>\mu}\lambda_i|\beta_i|^2\geq
         \mu^{1-s}\sum_{\lambda_i>\mu}\lambda_i^s|\beta_i|^2\geq\nonumber\\
          & \geq & K\mu^{1-s}
          \|\beta\|^2_{s,2} \geq
          K \tilde c^{-1}\mu^{1-s}|\beta|^2_p
        \end{eqnarray}
    and then
    by \eqref{non lin bound} and \eqref{control},
    \begin{eqnarray*}\label{geometric1}
        \widehat{J}(\beta) & = &
        |d\beta|_2^2-\int_Mf(\langle\beta+d\Phi(\beta),\beta+d\Phi(\beta)\rangle)\,\omega\\
        &\geq&|d\beta|_2^2-a |\beta|_p^p-b'|\beta|_p^2\\
        &\geq&\rho^2-a\left(\frac{\tilde
        c}{K\mu^{1-s}}\right)^\frac{p}{2}\rho^p-b'\frac{\tilde
        c}{K\mu^{1-s}}\rho^2\\
        &  = &\rho^2-a\left(\frac{\tilde c}{K}\right)^\frac{p}{2}-
        b'\frac{\tilde c}{K}.
    \end{eqnarray*}
    Since
        $$\lim_{\rho\rightarrow +\infty}K=1$$
    from the previous chain of inequalities we get \eqref{geometry1} for $\rho$
    large enough.

    Now take $\beta\in H^-(\mu)$. Observe that
        \begin{equation}\label{ineqonH-}
            |d\beta|_2^2\leq\lambda_k|
            \beta|_2^2.
        \end{equation}
    Moreover, since $\beta$ and $d\Phi(\beta)$ are orthogonal in the space $L^2_k(M)$,
    there exists $K_1>0$ s.t.
        \begin{eqnarray}\label{betortdPhibet}
            |\beta+d\Phi(\beta)|^p_p
            &  \geq  &
            K_1(|\beta+d\Phi(\beta)|_2^2)^\frac{p}{2}\nonumber\\
            &  \geq  &
            K_1(|\beta|_2^2+|d\Phi(\beta)|_2^2)^\frac{p}{2}\nonumber\\
            &  \geq  &
            K_1|\beta|_2^p.
        \end{eqnarray}
    By \eqref{increasing}, \eqref{ineqonH-} and \eqref{betortdPhibet},
    since in $H^-(\mu)$ all the norms are equivalent,
    there exists $ K_2,K_3>0$ such that
    \begin{eqnarray*}
        \widehat{J}(\beta) &=&
        |d\beta|_2^2-\int_Mf(\langle\beta+d\Phi(\beta),\beta+d\Phi(\beta)\rangle)\,\omega\\
        &\leq&\lambda_k|\beta|_2^2-c|\beta+d\Phi(\beta)|^p_p+d \mis(M) \\
        &\leq&\lambda_k|\beta|_2^2-K_2|\beta|_2^p+K_3\\
        &  \leq  &
        \sup_{t>0}(\lambda_kt^2-K_2t^p+K_3).
    \end{eqnarray*}
    \end{proof}
\smallskip
So we are ready to give the following
    \begin{proof}[Proof (of the first part of Theorem \ref{main}).]
        By Lemma \ref{ps condition}, Lemma
        \ref{geometry} and \eqref{3a}, using Theorem \ref{mult}, we
        find infinitely many couples of critical points. In fact, for all $i\geq1$
        there exist at least $\dim M_{\lambda_{k_i}}$
        critical points for $\widehat J$, whose critical values are in the interval
        $[C_{2i},C_{2i+1}].$ Since the sequence $(C_i)_{i\geq1}$
        is strictly increasing, we obtain a countable set of critical
        points.
    \end{proof}

\section{Zero mass case}\label{zeromass}
In this section we consider again the problem \eqref{problem},
replacing condition $f_1$ by $\tilde f_1$. \\
Moreover we replace $f_2$ by the technical hypothesis $\tilde f_2$
that, as already seen, implies for all $q\in M$ the strict
convexity of the functional
    $$I_q(\xi)=f(\langle\xi,\xi\rangle_q).$$
Observe that, integrating in $\tilde f_2$, by density we also
deduce that
    \begin{equation}\label{Fconv}
        F\hbox{ and } F_\beta\hbox{ are strictly convex on }  L_k^p(M).
    \end{equation}

Now, in order to prove the second part of Theorem \ref{main}, we
consider the perturbed equation
    \begin{equation}\label{pert}
        \delta d\xi=f'_\varepsilon(\langle\xi,\xi\rangle)\xi
    \end{equation}
where
    \begin{equation*}
        f_\varepsilon(t)=f(t)+\varepsilon t,\;\varepsilon>0.
    \end{equation*}
We set
    \begin{eqnarray}
        F_\varepsilon:
        &=&\xi\in L_k^p(M)\mapsto\int_Mf_\varepsilon(\langle\xi,\xi\rangle)\omega\\
        \noalign{\hbox{and for all $\beta\in
        L_k^p(M)$}}\nonumber\\
        (F_\varepsilon)_\beta: & = & \alpha\in\V\mapsto
        F_\varepsilon(d\alpha+\beta).
    \end{eqnarray}
The function $f_\varepsilon$ satisfies $f_1$, $f_3$ and $f_4$,
and, by \eqref{Fconv}, $F_\varepsilon$ and $(F_\varepsilon)_\beta$
are uniformly convex respectively on $L_k^p(M)$ and $\V$. From the
first part, we conclude that the equation \eqref{pert} possesses
infinitely many $\varepsilon-$solutions of the type
$\xi_\varepsilon=\beta_\varepsilon+d\Phi_\varepsilon(\beta_\varepsilon)$
where
\begin{equation}\label{Phivar}
\begin{array}{cll}
\Phi_\varepsilon: L^p_k(M)\rightarrow\V & s.t.
&\Phi_\varepsilon(\beta)\hbox{ is the minimizer of
}(F_\varepsilon)_\beta
\end{array}
\end{equation}
and $\beta_\varepsilon$ is a critical point of the functional
    \begin{eqnarray}\label{Jeps}
        \widehat J_\varepsilon(\beta)&=&\int_M\langle
        d\beta,d\beta\rangle\,\omega- \int_M f_\varepsilon(\langle
        d\Phi_\varepsilon(\beta)+\beta,d\Phi_\varepsilon(\beta)+\beta\rangle)\,\omega\nonumber\\
        &=&\int_M\langle
        d\beta,d\beta\rangle\,\omega- \int_M f(\langle
        d\Phi_\varepsilon(\beta)+\beta,d\Phi_\varepsilon(\beta)+
        \beta\rangle)\,\omega\nonumber\\
                                        & &\phantom{\int_M\langle
        d\beta,d\beta\rangle\,\omega}-\varepsilon\int_M\langle d\Phi_\varepsilon(\beta)
        +\beta,d\Phi_\varepsilon(\beta)+\beta\rangle\,\omega.
    \end{eqnarray}

Now we construct infinitely many sequences $(\beta_n)_n$ of
$\varepsilon_n-$solutions  of \eqref{pert}, in such a way, passing
to the limit, we get solutions for the non-perturbed
problem \eqref{problem}.\\
Of course, assuming that two different sequences converge, we have
no chance to prove that the corresponding limits are
different without any a-priori estimate on the critical values.\\
So we need a separation property on the sequences and, with
reference to this, we introduce the following
    \begin{defin}
        Let $(\varepsilon_n)_{n\geq1}$ be a sequence of real numbers and set
        $(\beta_n^1)_{n\geq1}$ and $(\beta_n^2)_{n\geq1}$ two
        sequences of $k$-forms. We say that $(\beta_n^1)_{n\geq1}$ and
        $(\beta_n^2)_{n\geq1}$ are well-separated if
        there exist $k_1$, $k_2$, $k_3$, $k_4\in\R$ such that
            $$k_1\leq \widehat J_{\varepsilon_n}(\beta_n^1)\leq
            k_2<k_3\leq\widehat J_{\varepsilon_n}(\beta_n^2)\leq
            k_4,\quad\forall n\geq 1.$$
    \end{defin}
Actually we have this result
    \begin{teo}\label{wellsep}
        Let $\varepsilon_n \searrow 0^+.$ Then there exists a
        countable set of sequences $\big((\beta^i_n)_{n\geq 1}\big)_{i\geq1}$ of
        $\varepsilon_n-$solutions
        that are each other well-separated.
    \end{teo}

    \begin{proof}
        Observe that, using Theorem \ref{mult} and its notations,
        we get the conclusion if we prove that
        there exist two sequences $(\mu_i)_{i\geq 0}$ and $(\rho_i)_{i\geq0}$, and a strictly
        increasing sequence $(C_i)_{i\geq 0}$ such that $3b$ and $3c$ of Lemma \ref{geometry}
        hold $\varepsilon_n-$uniformly, i.e.
        \begin{eqnarray}
            C_{2i}\leq \widehat J_{\varepsilon_n}(\beta),&&\forall
            \beta\in\partial S_{\rho_i}\cap
            H^+(\mu_i),\;\forall\varepsilon_n\\
            \widehat
            J_{\varepsilon_n}(\beta)<C_{2i+1},&&\forall\beta\in
            H^-(\mu_i),\;\forall\varepsilon_n.
        \end{eqnarray}\\
        We shall prove that for every
        $C>0$ there exist $\rho>0$ and $\mu>0$ such that
        \begin{eqnarray}
                \inf_{\substack{\beta\in\partial S_\rho\cap H^+(\mu)\\
                n\geq1}}\widehat J_{\varepsilon_n}(\beta)&\geq& C,\\
                \sup_{\substack{\beta\in H^-(\mu)\\n\geq1}}
                \widehat J_{\varepsilon_n}(\beta) & < &
                +\infty.
        \end{eqnarray}
        By $f_3$ and \eqref{Phivar}, using the embedding
        $L^p_k(M)\hookrightarrow L^2_k(M)$,
        and since $\varepsilon_n$ is decreasing we have
        \begin{equation}\label{estim1}
            \begin{array}{lcl}
                \varepsilon_n|\beta+d\Phi_{\varepsilon_n}(\beta)|_2^2 &
                + &\displaystyle\int_Mf(\langle\beta+d\Phi_{\varepsilon_n}(\beta),\beta+
                d\Phi_{\varepsilon_n}(\beta)\rangle)\,\omega\\
                &\leq&\varepsilon_n|\beta|_2^2+
                \displaystyle\int_Mf(\langle\beta,\beta\rangle)\,\omega\\
                &\leq&\varepsilon_1|\beta|_2^2+a
                |\beta|_p^p+b'|\beta|_p^2\\
                &\leq&a|\beta|_p^p+b''|\beta|_p^2\,,
            \end{array}
        \end{equation}
        where $b''$ is a suitable positive constant.

        Set $\mu=\rho^{\frac{2}{1-s}}$, and take $\beta\in H^+(\mu)\cap S_\rho$.
         By \eqref{estim1} and using \eqref{control},
        \begin{eqnarray}\label{firststep}
            \widehat J_{\varepsilon_n}(\beta) & = &
            |d\beta|_2^2-
            \varepsilon_n|\beta+d\Phi_{\varepsilon_n}(\beta)|_2^2-\nonumber\\
            &&\phantom{|d\beta|_2^2}-\int_Mf(\langle\beta+d\Phi_{\varepsilon_n}(\beta),\beta+
            d\Phi_{\varepsilon_n}(\beta)\rangle)\,\omega\nonumber\\
            &\geq&
            |d\beta|_2^2-a |\beta|_p^p-b''|\beta|_p^2\nonumber\\
            &\geq&\rho^2-a\left(\frac{\tilde c}{K}\right)^\frac{p}{2}-
            b''\frac{\tilde c}{K}\geq C
        \end{eqnarray}
    uniformly for $n\geq1$, for $\rho$ large enough.\\
    Moreover, if $\beta\in H^-(\mu)$, then by \eqref{increasing}
        \begin{equation}\label{estim2}
            \int_Mf(\langle\beta+d\Phi_{\varepsilon_n}(\beta),\beta+
            d\Phi_{\varepsilon_n}(\beta)\rangle)\,\omega\geq
            c|\beta+d\Phi_{\varepsilon_n}(\beta)|_p^p-d\mis(M).
        \end{equation}
        So, if we set $\lambda_{k}=\min\{\lambda_j\mid\lambda_j>\mu\},$
        then, for suitable $c_1,c_2>0$, by \eqref{estim2} we have

        \begin{eqnarray}
            \widehat J_{\varepsilon_n}(\beta)&=&|d\beta|_2^2-\varepsilon_n|\beta
            +d\Phi_{\varepsilon_n}(\beta)|_2^2-\nonumber\\
            &&\phantom{|d\beta|_2^2}-\int_Mf(\langle\beta+
            d\Phi_{\varepsilon_n}(\beta),\beta+d\Phi_{\varepsilon_n}(\beta)\rangle)\,\omega
            \leq\nonumber\\
            &\leq&\lambda_k|\beta|_2^2-c|\beta+d\Phi_{\varepsilon_n}(\beta)|_p^p+d\mis(M)
            \leq\nonumber\\
            &\leq&\lambda_{k}|\beta|_2^2-c_1(|\beta|_2^2+
            |d\Phi_{\varepsilon_n}(\beta)|_2^2)^{\frac{p}{2}}+c_2\leq\nonumber\\
            &\leq&\lambda_{k}|\beta|_2^2-c_1|\beta|_2^p+c_2\leq\nonumber\\
            &\leq&\sup_{t\geq
            0}(\lambda_{k}t^2-c_1t^p+c_2)<+\infty\nonumber
        \end{eqnarray}
        uniformly for $n\geq1.$
        \end{proof}

\subsection{Some preliminary results}\label{lemmas}
    \begin{lemma}\label{lemma1}
        Let $(\xi_n)_{n\geq1}$ a sequence of forms in $L_k^p(M)$
        and $\xi\in L^p_k(M).$\\
        If $f$ satisfies $\tilde f_2, f_3$ and
        \begin{eqnarray}
            \xi_n&\rightharpoonup&\xi\hbox{ in
            }L_k^p(M)\label{wconv}\\
            \int_Mf(\langle\xi_n,\xi_n\rangle)\,\omega&\rightarrow&
            \int_Mf(\langle\xi,\xi\rangle)\,\omega\label{convf}
        \end{eqnarray}
    then
        \begin{equation}\label{sconv}
            \xi_n\rightarrow\xi\hbox{ in }L_k^p(M).
        \end{equation}
    \end{lemma}
    \begin{proof}
        First suppose $\xi_n$ and $\xi$ in $\Lambda_k(M)$. By $\tilde f_2$ we have
        $$f(\langle \xi_n,\xi_n\rangle)-f(\langle \xi,\xi\rangle)-
            2 f'(\langle
            \xi,\xi\rangle)\langle\xi_n,\xi_n-\xi\rangle\geq\overline
            c\langle \xi_n-\xi,\xi_n-\xi\rangle^\frac{p}{2},$$
        so, integrating, we obtain,
        \begin{eqnarray}\label{near}
            \int_M\big (f(\langle\xi_n,\xi_n\rangle)-f(\langle\xi,\xi\rangle)\big )\,
            \omega-2\Psi_\xi(\xi_n-\xi)\geq
            \phantom{c_2\int_M|\langle\xi_n,\xi_n\rangle^{\frac{1}{2}}-\langle\xi,}
            \nonumber\\
            \geq\
            \overline c\int_M\langle\xi_n-\xi,\xi_n-\xi\rangle^{\frac{p}{2}}\,\omega
        \end{eqnarray}
        where $\Psi_\xi$ represents the map
        $$\Psi_\xi:\eta\in
        L^p_k(M)\mapsto\int_Mf'(\langle\xi,\xi\rangle)\langle\xi,\eta\rangle\,\omega.$$
        Observe that $\Psi_\xi$ is  linear and continuous by $f_3$ so we get
        \eqref{sconv} from \eqref{wconv}, \eqref{convf} and
        \eqref{near}.\\
        If $\xi_n$ and $\xi$ are in $L_k^p(M),$ then we get the
        same conclusion by density.
        \end{proof}
        \begin{lemma}\label{lemma2}
            $\forall \beta\in\W$, $\forall\varepsilon>0:$
        \begin{equation}\label{ineqlemma2}
            \begin{array}{lcl}
                0&\leq&\displaystyle\int_Mf(\langle\beta+d\Phi_\varepsilon(\beta),\beta+
                d\Phi_\varepsilon
                (\beta)\rangle)\,\omega-\\
                &&\phantom{\int_Mf\langle\beta+d\Phi_\varepsilon(\beta),\beta}
                -\displaystyle\int_Mf(\langle\beta+d\Phi(\beta),\beta+d\Phi(\beta)\rangle)
                \,\omega\leq\\
                &\leq&\varepsilon(|\beta+d\Phi(\beta)|_2^2-|\beta+
                d\Phi_\varepsilon(\beta)|_2^2)
            \end{array}
        \end{equation}
        \end{lemma}
        \begin{proof}
            Consider $\beta\in\W$ and $\varepsilon>0$, and set
            $b=\langle\beta+d\Phi(\beta),\beta+d\Phi(\beta)\rangle$
            and $b_\varepsilon=\langle\beta+d\Phi_\varepsilon(\beta),\beta+
            d\Phi_\varepsilon(\beta)\rangle$. By definition of
            $\Phi$ and $\Phi_\varepsilon$, we have that
            \begin{eqnarray}
                F(\beta+d\Phi(\beta))&\leq&
                F(\beta+d\Phi_\varepsilon(\beta))\\
                F_\varepsilon(\beta+d\Phi_\varepsilon(\beta))&\leq&
                F_\varepsilon(\beta+d\Phi(\beta))
            \end{eqnarray}
        that is
            \begin{equation}\label{bbeps}
                \int_Mf(b)\,\omega\leq\int_Mf(b_\varepsilon)\,\omega
            \end{equation}
            and
            \begin{equation}\label{bepsb}
                \int_Mf(b_\varepsilon)\,\omega+\varepsilon
                |\beta+d\Phi_\varepsilon(\beta)|_2^2
                \leq\int_Mf(b)\,\omega+\varepsilon|\beta+d\Phi(\beta)|_2^2.
            \end{equation}
            Combining \eqref{bbeps} and \eqref{bepsb} together, we
            get \eqref{ineqlemma2}.
        \end{proof}
    \subsection{Main Theorem (second part)}\label{proof2}
        The following lemma holds
        \begin{lemma}\label{lemma3}
            Let $\varepsilon_n\searrow 0^+$ and
            $(\beta_n)_{n\geq1}$ a sequence in $\W$ such that
            \begin{description}
                \item{a)} $\exists k_2>k_1>0$ s.t. $k_1\leq\widehat
                J_{\varepsilon_n}(\beta_n)\leq k_2,$ $\forall n\geq 1,$
                \item{b)} $\widehat J '_{\varepsilon_n}(\beta_n)=0,$
                $\forall n\geq 1.$
            \end{description}
            Then there exists $\beta\in\W$ and a subsequence relabelled $(\beta_n)_{n\geq1}$
             such that
            \begin{description}
                \item{i)}
                $\phantom{ii}d\Phi_{\varepsilon_n}(\beta_n)\longrightarrow d\Phi(\beta)$
                in $L_k^p(M)$
                \item{ii)} $\phantom{i}\beta_n\longrightarrow\beta$ in
                $H_k^1(M)$
                \item{iii)} $k_1\leq\widehat J(\beta)\leq k_2.$
            \end{description}
        \end{lemma}
        \begin{proof}
            To simplify the notations, set $b_n=\langle\beta_n+d\Phi_{\varepsilon_n}(\beta_n),
            \beta_n+
            d\Phi_{\varepsilon_n}(\beta_n)\rangle$.  By $a)$ we have
            \begin{equation}\label{lem3eq1}
                k_1\leq|d\beta_n|_2^2-\varepsilon_n|\beta_n+d\Phi_{\varepsilon_n}
                (\beta_n)|_2^2-
                \int_Mf(b_n)\,\omega\leq k_2
            \end{equation}
            and, on the other hand, by $b$) and using
            \eqref{partderivJ=0inPhibeta},
            \begin{equation}\label{lem3eq2}
                |d\beta_n|_2^2-\varepsilon_n|\beta_n+d\Phi_{\varepsilon_n}(\beta_n)|_2^2-
                \int_Mf'(b_n)b_n\,\omega=0.
            \end{equation}
            Using arguments similar to those in the proof of Theorem \ref{ps condition},
             we have that $(\beta_n)_{n\geq 1}$
            and $(\Phi_{\varepsilon_n}(\beta_n))_{n\geq 1}$ are
            bounded respectively in $\W$ and $\V$, so there exist
            $\beta\in\W$ and $\eta\in\V$ such that (up to a
            subsequence)
            \begin{eqnarray}
                \beta_n\rightharpoonup\beta&\hbox{in}&H_k^1(M)\label{lem3eqq2}\\
                \Phi_{\varepsilon_n}(\beta_n)\rightharpoonup\eta&\hbox{in}&H_{k-1}^{1,p}
                \label{lem3eq3}\\
                \noalign{\hbox{and by compactness}}
                \beta_n\rightarrow\beta&\hbox{in}&L_k^p(M).\label{lem3eq4}
            \end{eqnarray}
            Now applying \eqref{ineqlemma2} to $\beta_n$ for every
            $n\geq 1$, we obtain
            \begin{eqnarray}
                0&\leq&\int_Mf(b_n)\,\omega-\int_Mf(\langle\beta_n+d\Phi(\beta_n),\beta_n+
                d\Phi(\beta_n)
                \rangle\,\omega\leq\label{forgotteneq1}\\
                &\leq&\varepsilon_n(|\beta_n+d\Phi(\beta_n)|_2^2-|\beta_n+
                d\Phi_{\varepsilon_n}
                (\beta_n)|_2^2).\label{forgotteneq2}
            \end{eqnarray}
            We claim that
            \begin{equation}\label{lem3eq5}
                \varepsilon_n(|\beta_n+d\Phi(\beta_n)|_2^2-
                |\beta_n+d\Phi_{\varepsilon_n}(\beta_n)|_2^2)
                \stackrel{n\rightarrow+\infty}{\longrightarrow}0.
            \end{equation}
            In fact, suppose by contradiction that
            \eqref{lem3eq5} is not true. Since \eqref{lem3eq1} and
            \eqref{lem3eq2} imply that
            \begin{equation}\label{lem3eq6}
            |\beta_n+d\Phi_{\varepsilon_n}(\beta_n)|_2\hbox{ is
            bounded (see proof of Theorem \ref{ps condition})},
            \end{equation}
            by the contradiction hypothesis we have that, up to a
            subsequence,
                $$|\beta_n+d\Phi(\beta_n)|_2\rightarrow+\infty$$
            and then
            \begin{equation}\label{lem3eq7}
                |\beta_n+d\Phi(\beta_n)|_p\rightarrow+\infty.
            \end{equation}
            By \eqref{increasing} and
            \eqref{lem3eq7}
            \begin{equation}\label{lem3eq8}
                \int_Mf(\langle
                \beta_n+d\Phi(\beta_n),\beta_n+d\Phi(\beta_n)\rangle)\,\omega\rightarrow+
                \infty.
            \end{equation}
            Comparing \eqref{forgotteneq1} and \eqref{lem3eq8} we
            deduce that
            \begin{equation}\label{intfgoesinf}
                \int_Mf(b_n)\,\omega\rightarrow+\infty.
            \end{equation}
            On the other hand, from \eqref{lem3eq1} and \eqref{lem3eq2} we
            have that
            \begin{equation}\label{bound2}
                |d\Phi_{\varepsilon_n}(\beta_n)+\beta_n|_p\hbox{
                is bounded (see proof of Theorem
                \ref{ps condition})},
            \end{equation}
            so, considering $f_3$, we get
                \begin{equation*}
                    \int_Mf(b_n)\,\omega\hbox{ }\hbox{
                    is bounded}
                \end{equation*}
            that contradicts \eqref{intfgoesinf}.

            Now observe that  \eqref{forgotteneq1},
            \eqref{forgotteneq2} and \eqref{lem3eq5} imply that
            \begin{equation}\label{lem3eq9}
                \int_Mf(b_n)\,\omega-\int_Mf(\langle\beta_n+d\Phi(\beta_n),
                \beta_n+d\Phi(\beta_n)\rangle)\,\omega\longrightarrow0.
            \end{equation}
            Since the map
            \begin{equation*}
                \beta\in\W\mapsto\int_Mf(\langle\beta+d\Phi(\beta),\beta+
                d\Phi(\beta)\rangle)\,\omega
            \end{equation*}
            is weakly continuous (see \cite {B.F.} or \cite {A}),
            \eqref{lem3eqq2} implies that
            \begin{equation}\label{wcont}
                \int_Mf(\langle\beta_n+d\Phi(\beta_n),\beta_n+d\Phi(\beta_n)\rangle)\,\omega
                \rightarrow
                \int_Mf(\langle\beta+d\Phi(\beta),\beta+d\Phi(\beta)\rangle)\,\omega
            \end{equation}
                and then, by \eqref{lem3eq9},
            \begin{equation}\label{lem3eq10}
                \int_Mf(b_n)\,\omega\longrightarrow\int_Mf(\langle\beta+d\Phi(\beta),
                \beta+d\Phi(\beta)\rangle)\,\omega.
            \end{equation}
            Since $F$ is weakly lower semicontinuous, from
            \eqref{lem3eq3}, \eqref{lem3eq4} and \eqref{lem3eq10} we have
            \begin{eqnarray*}
                F_\beta(\eta)&=&F(\beta+d\eta)\\
                &\leq&\liminf_n
                F(\beta_n+d\Phi_{\varepsilon_n}(\beta_n))
                =F(\beta+d\Phi(\beta))=
                F_\beta(\Phi(\beta))
            \end{eqnarray*}
            and then, by the uniqueness of the minimizer of $F_\beta$,
            \begin{equation}\label{lem3eq14}
                \eta=\Phi(\beta).
            \end{equation}
            Now from \eqref{lem3eq3}, \eqref{lem3eq4} and
            \eqref{lem3eq14} we have that
            \begin{equation}\label{lem3eq15}
                \beta_n+d\Phi_{\varepsilon_n}(\beta_n)\rightharpoonup\beta+d\Phi(\beta)
                \hbox{ in }L_k^p(M)
            \end{equation}
            so, by \eqref{lem3eq10}, \eqref{lem3eq15} and Lemma
            \ref{lemma1}, we have that
            \begin{equation}\label{lem3eq16A}
                \beta_n+d\Phi_{\varepsilon_n}(\beta_n)\rightarrow
                \beta+d\Phi(\beta)\hbox{ in }L_k^p(M),
            \end{equation}
            and then, taking \eqref{lem3eq4} into account,
            \begin{equation}\label{lem3eq16}
                d\Phi_{\varepsilon_n}(\beta_n)\rightarrow
                d\Phi(\beta)\hbox{ in }L_k^p(M),
            \end{equation}
            that corresponds to the assertion $i).$

            Now we pass to the proof of $ii).$\\
            From $b)$ we have that
            \begin{equation}\label{lem3eq17}
                0=\widehat
                J'_{\varepsilon_n}(\beta_n)=L(\beta_n)-\varepsilon_n(\beta_n+d\Phi_
                {\varepsilon_n}
                (\beta_n))-K(\beta_n,d\Phi_{\varepsilon_n}
                (\beta_n))
            \end{equation}
            where $L$ is the Riesz isomorphism between $\W$
            and its dual and
            \begin{equation*}
                K:(\xi,\eta)\in H_k^1(M)\times L_k^p(M)\mapsto
                K(\xi,\eta)\in\big(H_k^1(M)\big)'.
            \end{equation*}
            From \eqref{lem3eq17} and considering \eqref{lem3eq6} we have
            \begin{equation*}
                L(\beta_n)-K(\beta_n,d\Phi_{\varepsilon_n}(\beta_n))\longrightarrow 0
            \end{equation*}
            and then
            \begin{equation}\label{lem3eq18}
                \beta_n-L^{-1}\big(K(\beta_n,d\Phi_{\varepsilon_n}(\beta_n))\big)
                \longrightarrow 0.
            \end{equation}
            Now observe that
            \begin{eqnarray}
                    K \hbox{ is compact with respect to } \xi\label{Kprop1}\\
                    K \hbox{ is continuous with respect to } \eta\label{Kprop2}
            \end{eqnarray}so, by \eqref{lem3eqq2} and
            \eqref{lem3eq16}, from \eqref{lem3eq18} we get (up to
            a subsequence)
            \begin{equation*}
                \beta_n\rightarrow L^{-1}\big(K(\beta,d\Phi(\beta))\big)
            \end{equation*}
            and hence $ii)$.

            Finally, note that from $ii)$, \eqref{lem3eq6} and
            \eqref{lem3eq10} we have
            \begin{equation}\label{lem3eq19}
                \widehat J_{\varepsilon_n}(\beta_n)\rightarrow
                \widehat J(\beta)
            \end{equation}
            so $iii)$ is a consequence of $a)$ and
            \eqref{lem3eq19}.
        \end{proof}

        And now we are ready for the following
        \begin{proof}[Proof (of the second part of Theorem \ref{main}).]
            Let $\varepsilon_n\searrow 0^+.$ By Theorem \ref{wellsep}, there exist
            infinitely many
            well separated sequences of the type described in $a)$
            and $b)$ of the Lemma \ref{lemma3}.\\
            Certainly each of these sequences (up to a subsequence)
            converges in $H_k^1(M)$ by $ii)$ and each limit is
            different from another by $iii).$\\
            Say $(\beta_n)_{n\geq1}$ one of these sequences and
            $\beta$ its limit. If we show that $\beta+d\Phi(\beta)$ is a
            solution for \eqref{weakproblem}, then we have
            finished.\\
            Let $\eta\in \Lambda^k(M).$ For every $n\geq1,$
            certainly
            \begin{equation}\label{lasteq1}
                \langle L(\beta_n),\eta\rangle=\varepsilon_n\int_M\langle\beta_n+
                d\Phi_{\varepsilon_n}
                (\beta_n),\eta\rangle\,\omega+\langle
                K(\beta_n,d\Phi_{\varepsilon_n}(\beta_n)),\eta\rangle
            \end{equation}
            where $L$ and $K$ are those defined in Lemma
            \ref{lemma3}.\\
            By continuity, $ii)$ of Lemma \ref{lemma3} implies
            that
            \begin{equation}\label{lasteq2}
                \langle L(\beta_n),\eta\rangle\longrightarrow\langle L(\beta),\eta\rangle
            \end{equation}
            while \eqref{Kprop1} and \eqref{Kprop2} together with
            $i)$ and $ii)$ of Lemma \ref{lemma3} imply
            \begin{equation}\label{lasteq3}
                \langle
                K(\beta_n,d\Phi_{\varepsilon_n}(\beta_n)),\eta\rangle\longrightarrow
                \langle
                K(\beta,d\Phi(\beta)),\eta\rangle.
            \end{equation}
            Since trivially
            \begin{equation*}
                \varepsilon_n\int_M\langle\beta_n+d\Phi_{\varepsilon_n}
                (\beta_n),\eta\rangle\,\omega\longrightarrow0,
            \end{equation*}
            by \eqref{lasteq1}, \eqref{lasteq2} and
            \eqref{lasteq3} we obtain
            \begin{equation}
                \int_M\langle
                d\beta,d\eta\rangle\,\omega=\int_Mf'(\langle\beta+d\Phi(\beta),
                \beta+d\Phi(\beta)\rangle)\langle\beta+d\Phi(\beta),\eta\rangle\,\omega
            \end{equation}
            and then the conclusion.
        \end{proof}

\section*{Appendix}
        In this appendix we want to show that assumption $\tilde
        f_2$ is satisfied by the function $f(t)=t^\frac{p}{2}$. We
        will prove it by the following more abstract result
        \begin{lem}
        Let $\Big(H,(\cdot|\cdot)\big)$ be an Hilbert space, and
        $\|\cdot\|$ the induced norm. If $p>2$, then there exists
        $\overline c>0$ s.t. for every $x,y\in H$ the following
        inequality holds
        \begin{equation}\label{conv cond}
            \|x\|^p-\|y\|^p-p\|y\|^{p-2}(y|x-y)\geq \overline c\|x-y\|^p.
        \end{equation}
        \begin{proof}
        In \cite{Kan} the following inequality has been proved
        for all $a,b\in\R$
        \begin{equation}\label{ineqforR}
            |a|^p-|b|^p-p|b|^{p-2}b(a-b)\geq K|a-b|^p,
        \end{equation}
        where $K>0$ does not depend on $a$ and $b$.
        Now we fix $y\in H$ and distinguish the following cases:
        \begin{itemize}
            \item $x=ty,$ $t\geq0$;
            \item $x=ty,$ $t<0$;
            \item $x\neq ty,$ $t\in\R$.
        \end{itemize}
        If $x=ty$ for $t\geq0$, then $(x|y)=\|x\|\|y\|$ and
        $(x-y|x-y)=(\|x\|-\|y\|)^2.$ So \eqref{conv cond} can be
        written as follows
        \begin{equation}
            \|x\|^p-\|y\|^p-p\|y\|^{p-1}(\|x\|-\|y\|)\geq\overline
            c\,\big|\|x\|-\|y\|\big|^p,
        \end{equation}
        that corresponds to \eqref{ineqforR} for $a=\|x\|$ and
        $b=\|y\|.$\\
        If $x=ty$ for $t<0$, then $(x|y)=-\|x\|\|y\|$ and
        $(x-y|x-y)=(\|x\|+\|y\|)^2.$ In this case, \eqref{conv cond} becomes
        \begin{equation}
            \|x\|^p-\|y\|^p+p\|y\|^{p-1}(\|x\|+\|y\|)\geq\overline
            c\,\big|\|x\|+\|y\|\big|^p,
        \end{equation}
        that corresponds to \eqref{ineqforR} for $a=\|x\|$ and
        $b=-\|y\|.$\\
        Finally, if $x\not\in\{ty|t\in\R\}$, then $x=x_1+x_2$
        where $x_1\in\{ty|t\in\R\}$ and $(x_2|y)=0$. Since \eqref{conv cond} holds for
        $x_1$, we have that there exist three positive constant $c_1$, $c_2$ and $c_3$ s.t.
        \begin{eqnarray*}
            \|x\|^p-\|y\|^p-p\|y\|^{p-2}(y|x-y)&=&(\|x_1\|^2+\|x_2\|^2)^\frac{p}{2}\\
            &&\phantom{(\|x_1\|^2}
            -\|y\|^p
            -p\|y\|^{p-2}(y|x_1-y)\\
            &\geq&\|x_1\|^p-\|y\|^p-p\|y\|^{p-2}(y|x_1)+\|x_2\|^p\\
            &\geq&c_1\big(\|x_1-y\|^p+\|x_2\|^p\big)\\
            &\geq&
            c_2\big(\|x_1-y\|^2+\|x_2\|^2\big)^\frac{p}{2}\\
            &=&c_3(\|x_1-y+x_2\|^2)^\frac{p}{2}=c_3\|x-y\|^p.
        \end{eqnarray*}
    \end{proof}
        \end{lem}

    Now, since $M$ is a compact Riemannian manifold, then for
    every $q\in M$ the space $\Lambda^k(T_q(M))$ of the $k-$forms
    at $q$ is an Hilbert space with the scalar product
    $\langle\cdot,\cdot\rangle_q$. So, by
    Lemma A.1, there exists $\overline c>0$ s.t. for every
    $\xi_q,\eta_q\in\Lambda^k(T_q(M))$ we have that
    \begin{equation}\label{proppointwise}
        \langle \xi_q,\xi_q\rangle_q^\frac{p}{2}-\langle \eta_q,\eta_q\rangle_q^\frac{p}{2}-
        p\langle
        \eta_q,\eta_q\rangle_q^{\frac{p}{2}-1}\langle\eta_q,\xi_q-\eta_q\rangle_q\geq
        \overline
        c\langle \xi_q-\eta_q,\xi_q-\eta_q\rangle_q^\frac{p}{2}.
    \end{equation}
    Since \eqref{proppointwise} holds pointwise, then for $\xi,\eta\in\Lambda^k(M)$
    the following inequality holds globally
    \begin{equation}
        \langle \xi,\xi\rangle^\frac{p}{2}-\langle \eta,\eta\rangle^\frac{p}{2}-
        p\langle
        \eta,\eta\rangle^{\frac{p}{2}-1}\langle\eta,\xi-\eta\rangle\geq\overline
        c\langle \xi-\eta,\xi-\eta\rangle^\frac{p}{2}.
    \end{equation}

        \vskip 1,5cm

\end{document}